\documentclass[letterpaper, 10 pt, conference]{ieeeconf} 

\IEEEoverridecommandlockouts                              
\usepackage{lipsum}
\usepackage{amsmath}
\usepackage{amssymb}
\usepackage{amsthm}
\usepackage{mathtools}
\usepackage{todonotes}
\usepackage{cite}

\usepackage{todonotes}

\theoremstyle{definition}
\newtheorem{definition}{Definition}

\newtheorem{assumption}{Assumption}
\newtheorem{problem}{Problem}

\usepackage{graphicx}
\graphicspath{ {resources/} }
\usepackage{booktabs}
\usepackage{float}
\usepackage{cite}
\let\labelindent\relax % fix enumitem package
\usepackage{enumitem}

\title{\LARGE \bf
	An Optimal Control Approach to Flocking
}

\author{Logan E. Beaver, \textit{Student Member, IEEE}, Chris Kroninger,\\ Andreas A. Malikopoulos, \textit{Senior Member, IEEE}
    \thanks{This research was supported by Combat Capabilities Development Command, Army Research Laboratory, MD, USA.}
	\thanks{Logan E. Beaver and Andreas A. Malikopoulos are with the Department of Mechanical Engineering, University of Delaware, Newark, DE 19711, USA.}%
	\thanks{Chris Kroninger is with Combat Capabilities Development Command, Army Research Laboratory, MD, USA. (emails: \tt\small{lebeaver@udel.edu};  \tt\small{christopher.m.kroninger.civ@mail.mil};  \tt\small{andreas@udel.edu}.)}%
}

\begin{document}
	
	\maketitle
	
	\begin{abstract}
	    Flocking behavior has attracted considerable attention in multi-agent systems. The structure of flocking has been predominantly studied through the application of artificial potential fields coupled with velocity consensus.
	    These approaches, however, do not consider the energy cost of the agents during flocking, which is especially important in large-scale robot swarms.
	    This paper introduces an optimal control framework to induce flocking in a group of agents.
	    Guarantees of energy minimization and safety are provided, along with a decentralized algorithm that satisfies the  optimality conditions and can be realized in real time.
	    The efficacy of the proposed control algorithm is evaluated through simulation in both MATLAB and Gazebo.
	\end{abstract}
	
	\section{Introduction}
	
    \subsection{Background}

    Complex systems consist of many interdependent, independent agents that are connected to each other in both space and time \cite{Malikopoulos2016c}. The individual interactions between agents may result in unpredictable emergent behavior at a large scale. As the world becomes increasingly complex \cite{Malikopoulos2015}, modern control approaches will be required to optimize individual agents with an eye on overall system behavior \cite{Malikopoulos, Malikopoulos2015b}.
    
    Multi-agent systems have attracted considerable attention in many applications due to their natural parallelization, general adaptability, and ability to self-organize \cite{Oh2017}.
    This has proven useful in many applications, such as transportation \cite{Malikopoulos2018}, construction \cite{Lindsey2012ConstructionTeams}, and surveillance \cite{Corts2009}.
    Controlling emergent flocking behavior has been of particular interest to robotics researchers since the seminal paper by Reynolds \cite{Reynolds1987}, which introduced three heuristic rules for flocking in digital animation: move toward neighboring flockmates, avoid collisions, and match velocity with neighbors.
    Flocking has many practical applications, such as mobile sensing networks, coordinated delivery, reconnaissance, and surveillance  \cite{Olfati-Saber2006FlockingTheory}.

%    However, fielding large robot swarms imposes significant cost constraints on each individual. These constraints limit each agent's computational and sensing ability and reduce their individual energy storage capacity.
%    This is the driving force behind energy-optimal control algorithms, which are necessary to extend the useful life of each agent. 

    Currently, the most popular approach to impose flocking in multi-agent systems is the use of artificial potential fields and velocity consensus \cite{Barve, Chung2018ARobotics, Oh2015}.
    In these approaches, flock aggregation and collision avoidance are both handled by the potential field, the design of which is still an open question \cite{Vasarhelyi2018OptimizedEnvironments}.
    Several theoretical guarantees have been proven for these types of flocking models \cite{Tanner2007}; however, minimum-energy flocking still remains relatively unexplored.
    
    In this paper, we propose an extension of our previous work on energy-optimal trajectories for formations \cite{Beaver2019AGeneration} and apply it to the problem of flocking. We provide the following three contributions:
	(1) a continuous decentralized optimal control framework for minimum-energy flocking, (2) a closed-form solution for the energy-minimizing control input in the centralized case, and (3) an adaptation of the centralized optimal solution to the decentralized case.
    
    There have been several approaches in the literature that have considered optimal flocking using dynamic programming for aircraft aggregation \cite{Quintero2013FlockingApproach}, and optimal control over discrete time \cite{Zhang2015ModelConstraints}. In contrast, our approach is continuous, has an analytical solution, and allows flexibility in the flock shape.
	Gomez et al. \cite{Gomez2016Real-TimeSwarms} used a centralized controller to derive optimal trajectories for teams of agents to flock. In contrast, our approach adopts the centralized solution to a decentralized system without requiring a central computer to plan trajectories.
	The framework we present in this paper is related to the robot ecology paradigm for long-duration autonomy \cite{Egerstedt2018RobotAutonomy,Ibuki2020Optimization-BasedBodies}. However, we apply optimal control over a planning horizon, rather than reacting to the environment with gradient flow.

%	Finally, Hu et al. \cite{Hu2018Self-triggeredSystems} introduced a discretized model predictive control framework to minimize energy consumption, avoid collisions, and reduce communication via self-triggered updates.
%	Our work provides an analytical solution to the continuous form of this problem, reducing the computational burden on each agent and allowing our algorithm to be executed in real-time.

%	The most significant drawback of potential field methods is that the resulting trajectories are sub-optimal from an energy perspective.
%	This was explicitly stated by Reynolds, who observed that potential fields tend to apply forces perpendicular to flock and agent motion \cite{Reynolds1987}.
%	Several modern approaches have proposed alternative strategies to avoid inter-agent collisions, such as reactive obstacle avoidance \cite{Egerstedt2001FormationControl}, dynamic state constraints \cite{Pereira, Luis2019TrajectoryControl, Morgan2016}, initial system constraints \cite{Turpin2014}, and ordered/sequential trajectory planning \cite{Turpin}.
%	Work has been completed with respect to pairwise optimization constraints \cite{Bhattacharya2018DistributedPlanning}, but it has not been applied to flocking.

%	\subsection{Contributions of This Paper}

	The structure of the paper is as follows. We formulate the problem in Section \ref{sc:formulation}. In Section \ref{sc:solution}, we solve the optimization problem for the unconstrained and constrained cases. In Section \ref{sc:simulation}, we provide two sets of simulation results and discuss the observed emergent behavior. In the first set, we use MATLAB to demonstrate the viability of the proposed control scheme; in the second set, we use Gazebo to validate a minimum-communication approach, which reformulates the control method to use pure sensing. Finally, we draw conclusions and discuss future work in Section \ref{sc:conclusion}.
	
	\section{Problem Formulation} \label{sc:formulation}
	
	Consider a flock of  $N\in\mathbb{N}$ agents indexed by the set $\mathcal{A} = \{ 1, 2, \dots, N\}$. Each agent $i\in\mathcal{A}$ follows the double integrator dynamics,
	\begin{align}
	\dot{\mathbf{p}}_i(t) &= \mathbf{v}_i(t), \label{eq:pDynamics} \\
	\dot{\mathbf{v}}_i(t) &= \mathbf{u}_i(t), \label{eq:vDynamics}
	\end{align}
	where $t\in\mathbb{R}_{\geq0}$ is the time, and $\mathbf{p}_i(t), \mathbf{v}_i(t), \mathbf{u}_i(t) \in \mathbb{R}^2$ are the position, velocity, and control input, respectively. Each agent occupies a closed disk of radius $R\in\mathbb{R}_{>0}$. The state of each agent is given by
	\begin{equation} \label{eq:state}
	\mathbf{x}_i(t) = 
	\begin{bmatrix}
	\mathbf{p}_i(t) \\
	\mathbf{v}_i(t)
	\end{bmatrix}.
	\end{equation}
	The speed and control input are constrained such that
	\begin{align}
	||\mathbf{v}_i(t)|| \leq v_i^\text{max}, \label{eq:vConstraint}\\
	||\mathbf{u}_i(t)|| \leq u_i^\text{max}, \label{eq:uConstraint}
	\end{align}
	for all $t\in\mathbb{R}_{\geq0}$.
	
\begin{figure}[ht]
%	\vspace{0.5em}
	    \centering
	    \includegraphics[width=0.4\textwidth]{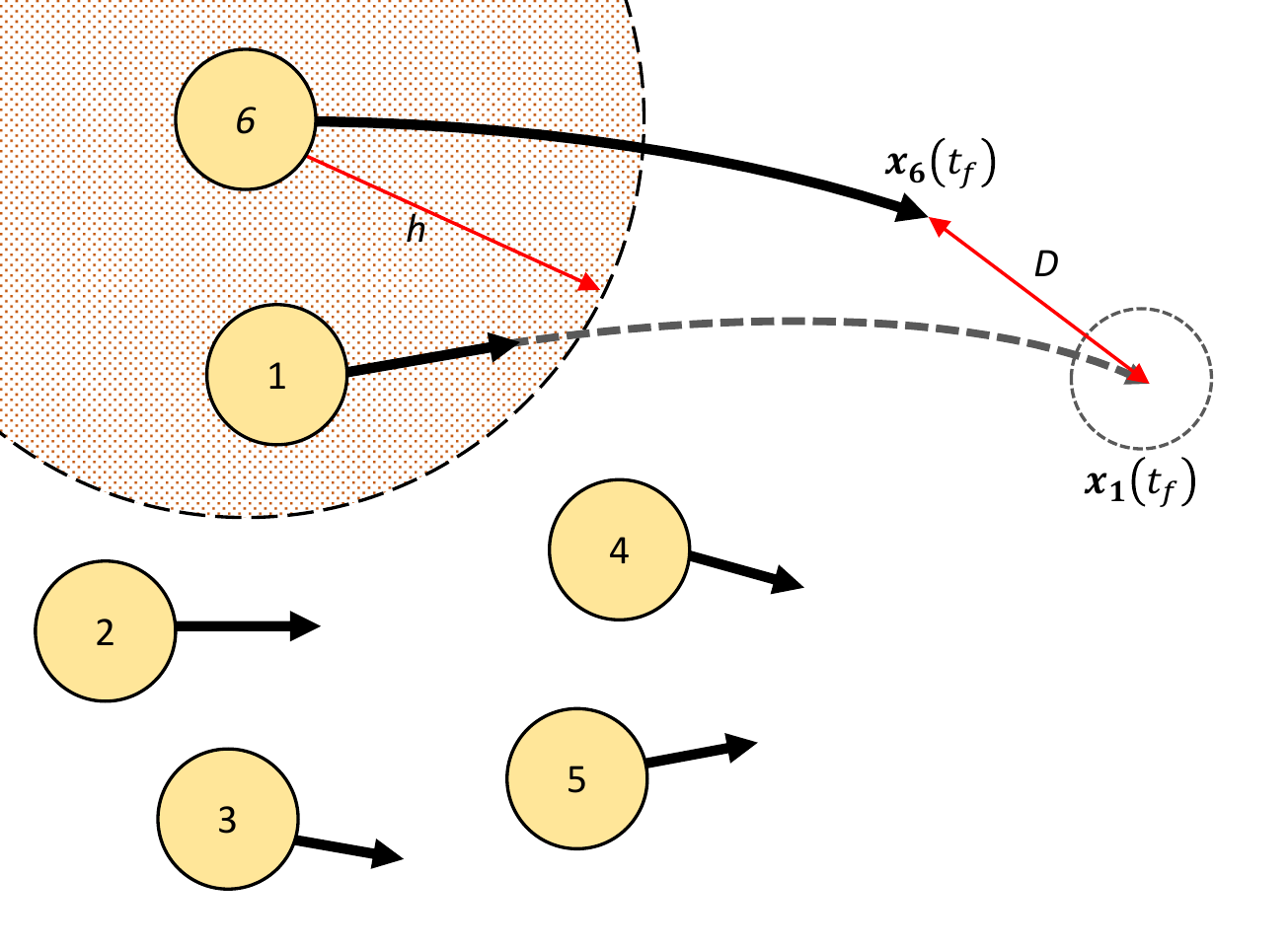}
	    \caption{A diagram of a flocking consisting of six agents.}
	    \label{fig:flocking}
	\end{figure}

	For any pair of agents $i,j\in\mathcal{A}$, the relative displacement between them is described by the vector
	\begin{equation} \label{eq:s}
	\mathbf{s}_{ij}(t) = \mathbf{p}_j(t) - \mathbf{p}_i(t), ~~ i,j\in\mathcal{A}.
	\end{equation} 
	To guarantee safety within the system, we impose the following constraints:
	\begin{align}
		\mathbf{s}_{ij}(t)\cdot\mathbf{s}_{ij}(t) &\geq 4R^2, ~~ \forall\, j\in\mathcal{A}, ~~ \forall t\in\mathbb{R}_{\geq0},\label{eq:noCollide} \\
		h &> 2R, \label{eq:sensingConstraint}
	\end{align}
	where \eqref{eq:noCollide} guarantees collision avoidance and \eqref{eq:sensingConstraint} is a system-level constraint which allows collisions to be detected before they occur. %The squared form of \eqref{eq:noCollide} is selected as its derivative with respect to $\mathbf{p}_{ij}(t)$ is smooth, which will simplify the analytical solution presented in Section \ref{sc:solution}.
	Each agent has also a sensing/communicating distance, $h\in\mathbb{R}_{>0}$, which is used to define its neighborhood.
	\begin{definition} \label{df:neighborhood}
	The \emph{neighborhood} of each agent $i\in\mathcal{A}$, is defined by the set
	\begin{equation} \label{eq:neighborhood}
	\mathcal{N}_i(t) = \{ j \in\mathcal{A} ~|~  ||\mathbf{s}_{ij}(t)|| < h \},
	\end{equation}
	where $||\cdot||$ is the Euclidean norm.
	\end{definition}
	The neighborhood is allowed to switch over time, while $i\in\mathcal{N}_i(t)$ always holds. Agent $i$ is able to communicate with any agent $j\in\mathcal{N}_i(t)$ and sense its current state, $\mathbf{x}_{j}(t)$. A schematic of the system is presented in Fig. \ref{fig:flocking}.
	
	Finally, every agent is also equipped with a long-range sensor, which can estimate the centroid of the flock, given by
	\begin{equation} \label{eq:centroid}
	\mathbf{p}_{cg}(t) = \frac{1}{|\mathcal{N}_i(t)|} \sum_{i\in\mathcal{A}} \mathbf{p}_i(t),
	\end{equation}
	where $|\cdot|$ denotes set cardinality.
	The purpose of \eqref{eq:centroid} is to drive isolated agents back toward the flock. This can reasonably be achieved with inexpensive audio or visual sensors in the case where a majority of the agents have already formed an aggregate. However, for the case where there is no clear flock, it may be challenging to determine \eqref{eq:centroid} in practice. 
	
	In order to guarantee safety \eqref{eq:noCollide}, and to calculate \eqref{eq:centroid}, agent $i\in\mathcal{A}$ requires information about the trajectories of all agents $j\in\mathcal{N}_i(t).$ One potential approach to handle this is to impose a priority ordering on the agents; in this case, higher priority agents act first while lower priority agents must steer to avoid them \cite{Beaver2020AnAgents}. Our approach is for each agent to apply model predictive control. Each agent generates an initial unconstrained trajectory for a given time interval $t\in[t_0, t_f], $ where $t_f - t_0 > 0$ is the horizon. Then, agent $i$ follows this generated trajectory for some period $\Delta T$, at which point the trajectory is recalculated over the new interval $t\in[t_0 + \Delta T, t_f + \Delta T]$.
	
	In this case, a longer time interval corresponds to a larger planning space at the cost of increased computational complexity. A decrease in the replanning period, $\Delta T$,  provides agent $i$ with a better estimate of its neighbors' trajectories at the cost of more frequent calculations. The trajectory generated by each agent is the outcome of minimizing the flocking error function, defined next.
	\begin{definition} \label{df:flockError}
		For all agents $i\in\mathcal{A}$, the \emph{flocking error} is defined by the scalar function
		\begin{align} \label{eq:flockError}
		\Phi(\mathbf{x}_i, t_f) &=  w_1~\phi_d(\mathbf{x}_i, t_f) +
		w_2~\phi_v(\mathbf{x}_i, t_f) \notag\\ &\hspace{1em}+ w_3~\phi_a(\mathbf{x}_i, t_f), \\
		\phi_d(\mathbf{x}_i, t_f) &= ||\mathbf{v}_i(t_f) - \mathbf{v}_d(t_f) ||^2,
		\label{eq:velocityControl}\\
		\phi_v(\mathbf{x}_i, t_f) &= ||\mathbf{v}_i(t_f) - \mathbf{v}_{\text{avg}}(t_f)||^2,
		\label{eq:velocityMatching}\\
		\phi_a(\mathbf{x}_i, t_f) &= 
		\begin{cases}
		(||\mathbf{p}_{cg}(t_f) - \mathbf{p}_i(t_f)|| - D)^2, \hspace{0.1em} |\mathcal{N}_i(t_0)| = 1,\\
		\underset{j\in\mathcal{N}_i(t_0)}{\sum} \Big( \mathbf{s}_{ij}(t_f)\cdot\hat{p}_{ij}(t_f)-D  \Big)^2, \\
		\hspace{12em} |\mathcal{N}_i(t_0)| > 1,
		\end{cases}
		\label{eq:aggregation}
		\end{align}
		where $w_1, w_2, $ and $w_3$ are normalizing parameters which weight the influence of velocity control $\phi_d(\mathbf{x}_i, t_f)$, velocity matching $\phi_v(\mathbf{x}_i, t_f),$ and aggregation $\phi_a(\mathbf{x}_i, t_f)$ in the overall system behavior. The system parameters $\mathbf{v}_d(t_f)$ and $D$ set the desired flock velocity and separating distance between agents. The average agent speed, $\mathbf{v}_{\text{avg}}(t_f)$, is given by
		\begin{equation}
		\mathbf{v}_{\text{avg}}(t_f) = \begin{cases}
		\dot{\mathbf{p}}_{cg}(t_f),   & |\mathcal{N}_i(t_0)| = 1,\\
		\frac{1}{|\mathcal{N}_i(t_0)|} \sum_{j\in\mathcal{N}_i(t_0)} \mathbf{v}_j(t_f),  & |\mathcal{N}_i(t_0)| > 1.
		\end{cases} \label{eq:vAvg}
		\end{equation}
	\end{definition}
	
	The cases in \eqref{eq:aggregation} and \eqref{eq:vAvg} when $|\mathcal{N}_i(t_0)| = 1$ are equivalent to placing a virtual agent at the centroid of the flock for any agent which becomes isolated. This will drive isolated agents toward the centroid by the construction of the flocking error. However, this does not guarantee that a disconnected group of agents will return to the flock if separation occurs.
	
	Each agent generates its energy-optimal trajectory by solving the following decentralized optimal control problem.
	
	\begin{problem}{(\emph{Minimum-energy flocking})}\label{pr:energyMin}
	
	For every agent $i\in\mathcal{A}$,
	\begin{align}
	\min_{\mathbf{u}_i(t)} &\Big\{ \Phi(\mathbf{x}_i, t_f) + \int_{t_0}^{t_f} ||\mathbf{u}_i(t)||^2 dt \Big\}\\
	\text{subject to:} \hspace{1em} & \eqref{eq:pDynamics}, \eqref{eq:vDynamics}, \eqref{eq:vConstraint}, \eqref{eq:uConstraint}, \eqref{eq:noCollide}, ~ \mathbf{x}_i(t_0) = \mathbf{x}_i^0,\notag
	\end{align}
	where $\Phi(\mathbf{x}_i, t_f)$ is given by Definition \ref{df:flockError} and $\mathbf{x}_i^0$ is the initial state of agent $i$.
	\end{problem}
	By minimizing the $L^2$ norm of the control input we expect to see a proportional reduction in energy consumption.
	To solve Problem \ref{pr:energyMin}, we impose the following assumptions.
	\begin{assumption}\label{as:empty}
        There are no external disturbances or obstacles.
    \end{assumption}
    Assumption \ref{as:empty} is imposed to evaluate the idealized performance of the proposed algorithm. %It is well known that optimal control can be sensitive with respect to noise and disturbances, and 
    This assumption may be relaxed by introducing a measure of robustness into Problem \ref{pr:energyMin}.% Obstacle avoidance may also be added as an additional safety constraint on each agent to navigate around obstacles.
    
	\begin{assumption}\label{as:perfect}
	    There are no errors or delays with respect to communication and sensing.
	\end{assumption}
    The strength of Assumption \ref{as:perfect} is application dependent. In general, it has been shown that sparse updates to model predictive control may be sufficient \cite{Hu2018Self-triggeredSystems}. %However, these delays may become significant for large or fast-moving swarms in constrained environments.
    
    \begin{assumption}\label{as:lowDensity}
        The flock is low-density (only two agents $i,j\in\mathcal{A}$ ever come within distance $|\mathbf{s}_{ij}| = 2R$ of each other).
    \end{assumption}
	Assumption \ref{as:lowDensity} may be strong, but it is imposed to simplify the solution to the optimal speed profile for two agents that are safety constrained.
	This assumption may be removed if an order is imposed on the agents rather than using model predictive control \cite{Beaver2020AnAgents}, and it may be relaxed in cases when a numerical solver can generate the safety-constrained trajectory in real-time.

	\section{Solution Approach} \label{sc:solution}

    To derive an analytical solution for Problem \ref{pr:energyMin} we use Hamiltonian Analysis.
    As the case where the state and control constraints, \eqref{eq:vConstraint} and \eqref{eq:uConstraint}, are active is well studied in the literature \cite{Malikopoulos2018} we will only consider the safety constrained case.
	
	As a first step, for an agent $i\in\mathcal{A}$, the safety constraint, \eqref{eq:noCollide}, must be derived until the control input, $\mathbf{u}_i(t)$, appears,
	\begin{equation} \label{eq:tangency}
	\mathbf{N}_i(t) = \begin{bmatrix}
	4R^2 - \mathbf{s}_{ij}(t)\cdot\mathbf{s}_{ij}(t)\\
	-\mathbf{s}_{ij}(t)\cdot\dot{\mathbf{s}}_{ij}(t)\\
	-\mathbf{s}_{ij}(t)\cdot\ddot{\mathbf{s}}_{ij} - \dot{\mathbf{s}}_{ij}(t)
	\cdot\dot{\mathbf{s}}_{ij}(t)
	\end{bmatrix} \leq \mathbf{0}.
	\end{equation}
    The Hamiltonian is then augmented by the final row of $\mathbf{N}_i(t)$, which yields
	\begin{align} 
	H_i = ||\mathbf{u}_i(t)||^2 + \boldsymbol{\lambda}_i^p(t)\cdot\mathbf{v}_i(t) + \boldsymbol{\lambda}_i^v(t)\cdot\mathbf{u}_i(t) \notag \\
	- \sum_{j\in\mathcal{N}_i} \mu_{i,j}(t) \Big( \mathbf{s}_{ij}(t) \cdot \ddot{\mathbf{s}}_{ij}(t) + \dot{\mathbf{s}}_{ij}(t)\cdot\dot{\mathbf{s}}_{ij}(t)  \Big), \label{eq:hamiltonian}
	\end{align}
	where $\boldsymbol{\lambda}_i^p(t), \boldsymbol{\lambda}_i^v(t)$ are the position and velocity covectors, and $\mu_{ij}(t)$ is the Lagrange multiplier with values
	\begin{align}
		\mu_{ij}(t) = 
		\begin{cases}
			\geq 0,  & \text{if  } {\mathbf{s}_{ij}(t) \cdot \ddot{\mathbf{s}}_{ij}(t) + \dot{\mathbf{s}}_{ij}(t)\cdot\dot{\mathbf{s}}_{ij}(t)} = 0, \\
			0,  & \text{if  } {\mathbf{s}_{ij}(t) \cdot \ddot{\mathbf{s}}_{ij}(t) + \dot{\mathbf{s}}_{ij}(t)\cdot\dot{\mathbf{s}}_{ij}(t) > 0}.
		\end{cases}
	\end{align}

	To solve \eqref{eq:hamiltonian} for agent $i\in\mathcal{A}$ we consider that (i) all agents $j\in\mathcal{N}_i(t)$ satisfy $\mu_{i,j} = 0$ or 
	(ii) any agent $j\in\mathcal{N}_i(t)$ satisfies $\mu_{i,j} > 0$. We then piece the constrained and unconstrained arcs together to arrive at a piecewise-continuous, energy-optimal trajectory. Next, we present the unconstrained motion and boundary conditions followed by our algorithm for constructing an energy-optimal trajectory in real-time.	
	
	\subsection{Unconstrained Motion} \label{ss:unconstrained}
	The energy-optimal unconstrained trajectories for the position, speed, acceleration, and covectors resulting from \eqref{eq:hamiltonian} over an interval $t\in[t_1, t_2]\subset\mathbb{R}_{\geq0}$ are \cite{Malikopoulos2018}
	\begin{align}
		\mathbf{p}_i(t) =&~ \frac{1}{6}\mathbf{a}_i t^3 + \frac{1}{2}\mathbf{b}_i t + \mathbf{c}_i t + \mathbf{d}_i,\label{eq:uncP}\\
		\mathbf{v}_i(t) =&~ \frac{1}{2}\mathbf{a}_i t^2 + \mathbf{b}_i t + \mathbf{c}_i \label{eq:uncV},\\
	    \mathbf{u}_i(t) =&~ \mathbf{a}_i t + \mathbf{b}_i \label{eq:uncU},\\
	    \boldsymbol{\lambda}^p_i(t) =&~ \mathbf{a}_i \label{eq:uncLp}, \\
	    \boldsymbol{\lambda}^v_i(t) =&~ -\mathbf{a}_i t - \mathbf{b}_i, \label{eq:uncLv}
	\end{align}
	where $\mathbf{a}_i, \mathbf{b}_i, \mathbf{c}_i, \mathbf{d}_i \in\mathbb{R}^2$ are constants of integration. Equations \eqref{eq:uncP} - \eqref{eq:uncLv} contain $8$ unknowns, which are solved using $8$ boundary conditions;
	\begin{align} 
	    \mathbf{x}_i(t_0) &= \mathbf{x}_i^0, \label{eq:IC} \\
	    \boldsymbol{\lambda}_i(t_f) &= \frac{\partial \Phi(\mathbf{x}_i,t)}{\partial \mathbf{x}_i}\Big|_{t_f}. \label{eq:phiBoundary}
	\end{align}
	
	Solving \eqref{eq:phiBoundary} yields
	\begin{align}
	\boldsymbol{\lambda}_i^p(t_f) &=
	-2w_3 \sum_{j\in\mathcal{N}_i(t_0)} \Big( \big( ||\mathbf{s}_{ij}(t_f)||-D \big) \,\mathbf{s}_{ij}(t_f) \Big), \label{eq:bc3}
	\\
	\boldsymbol{\lambda}_i^v(t_f) &= 2w_2\Big(\mathbf{v}_i(t_f)-\mathbf{v}_{\text{avg}}\Big)  + 2w_1\Big(\mathbf{v}_i(t_f)-\mathbf{v}_d\Big). \label{eq:bc4}
	\end{align}
	We then substitute \eqref{eq:s}, \eqref{eq:uncLp}, and \eqref{eq:uncLv} into \eqref{eq:bc3} and \eqref{eq:bc4} to solve for $\mathbf{p}_i(t_f)$ and $\mathbf{u}_i(t_f)$. This yields two equations
	\begin{align}
	\mathbf{p}_i(t_f) &= \frac{\mathbf{a}_i}{2w_3\sum_{j}(||\mathbf{s}_{ij}(t_f)|| - D)} + \sum_{j\in\mathcal{N}_i}\mathbf{p}_{j}(t_f),  \label{eq:boundary3}\\
	\mathbf{u}_i(t_f) &= -2w_2\Big(\mathbf{v}_i(t_f)-\mathbf{v}_{\text{avg}}\Big)  - 2w_1\Big(\mathbf{v}_i(t_f)-\mathbf{v}_d\Big), \label{eq:boundary4}
	\end{align}
	where $\mathbf{a}_i = 0$ if the right-hand side of \eqref{eq:bc3} is ever zero.
	
	Equations \eqref{eq:IC}, \eqref{eq:phiBoundary}, and \eqref{eq:boundary4} give three conditions to solve for the constants $\mathbf{b}_i$, $\mathbf{c}_i$, and $\mathbf{d}_i$ in \eqref{eq:uncP} - \eqref{eq:uncLv}. The value of $\mathbf{a}_i$ can then be found by the solution of \eqref{eq:boundary3}, which must be computed numerically.

	\subsection{Constrained Motion}
	To generate the safety-constrained motion of agent $i\in\mathcal{A}$, we require agent $i$ to cooperate with all agents $j\in\mathcal{N}_i(t)$ to solve the centralized optimal control problem whenever the collision avoidance constraint becomes active. The agents will then all employ the centralized solution to guarantee collision avoidance. %This is equivalent to minimizing the sum of the individual Hamiltonians, \eqref{eq:hamiltonian}, as the objectives of each agent are independent.
	
	For any agent $i\in\mathcal{A}$, we define the set $\mathcal{V}_i$ as
	\begin{equation} \label{eq:activeSet}
	\mathcal{V}_i(t) \coloneqq \{j\in\mathcal{N}_i(t) ~|~ \mu_{ij}(t) > 0, ~ j\neq i  \}.
	\end{equation}
	Thus, when $\mathcal{V}_i \neq \emptyset$ agent $i$ must follow a safety-constrained trajectory. Application of the Euler-Lagrange equations yields
	\begin{align}
	\mathbf{u}_i(t) =&~ -\boldsymbol{\lambda}_i^v(t) - \sum \mu_{ij}(t)\mathbf{s}_{ij}(t), \label{eq:uel} \\
	-\dot{\boldsymbol{\lambda}}_i^v(t) =&~ \boldsymbol{\lambda}_i^p(t) + \sum \mu_{ij}(t)\dot{\mathbf{s}}_{ij}(t), \label{eq:lambdaVdot}\\
	-\dot{\boldsymbol{\lambda}}_i^p(t) =&~ \sum \mu_{ij}(t)\ddot{\mathbf{s}}_{ij}(t), \label{eq:lambdaPdot}%\\
	%\mu_{ij}(t) =&~ \Big(\boldsymbol{\lambda}_j(t) - \boldsymbol{\lambda}_i(t)\Big)~\mathbf{s}_{ij}(t) - ||\dot{\mathbf{s}}_{ij}(t)||^2 \notag\\ &+ \sum_{k\in\mathcal{V}_j(t), k\neq i} \mu_{jk}(t)~\mathbf{s}_{jk}(t) \notag\\ &+ \sum_{k\in\mathcal{V}_i(t), k\neq j} \mu_{ik}(t)~\mathbf{s}_{ik}(t) \label{eq:muij} 
	\end{align}
%	\todoLogan{The equation for muij is incorrect, but the incorrect parts are assumed away to zero.}
	which, in general, must be solved numerically. However, under Assumption \ref{as:lowDensity}, we may consider the case where only two agents interact. 
	We define
	\begin{equation}
	    a_{ij}(t) \coloneqq ||\dot{\mathbf{s}}_{ij}(t)||, \label{eq:a}
	\end{equation}
	which is the relative speed between the two constrained agents. 
	We may then construct a new basis for $\mathbb{R}^2$ which we define next.
	\begin{definition}\label{df:contactBasis}
	The orthonormal \emph{contact basis} for any two agents in contact, $i\in\mathcal{A}, ~ j \in\mathcal{V}_i$, with a nonzero relative speed, is defined as
	\begin{align}
	\hat{p}_{ij}(t) \coloneqq&~ \frac{\mathbf{s}_{ij}(t)}{||\mathbf{s}_{ij}(t)||} = \frac{\mathbf{s}_{ij}(t)}{2R} \label{eq:phat},\\
	\hat{q}_{ij}(t) \coloneqq&~ \frac{\dot{\mathbf{s}}_{ij}(t)}{||\dot{\mathbf{s}}_{ij}(t)||} = \frac{\dot{\mathbf{s}}_{ij}(t)}{a_{ij}(t)}, \label{eq:qhat}
	\end{align}
	where $\hat{p}_{ij}(t)\cdot\hat{q}_{ij}(t) = 0$ by \eqref{eq:tangency}.
	\end{definition}
	
	%--- preliminary math
    To solve the \eqref{eq:uel} - \eqref{eq:lambdaPdot} we will project $\ddot{\mathbf{s}}_{ij}(t)$ onto the contact basis (Definition \ref{df:contactBasis}). From \eqref{eq:tangency} we can project $\ddot{\mathbf{s}}_{ij}$ onto $\hat{p}_{ij}$ by
	\begin{align}
	    \ddot{\mathbf{s}}_{ij}(t)\cdot\mathbf{s}_{ij}(t) = -\dot{s}_{ij}(t)\cdot\dot{s}_{ij}(t) = -a_{ij}^2(t). \label{eq:sddots}
	\end{align}
%	Next, we will solve for the projection of $\ddot{\mathbf{s}}_{ij}$ onto $\hat{q}_{ij}$ as follows: first, integrate this quantity using integration by parts
%    \begin{align}
%    \int \ddot{\mathbf{s}}_{i,j}(t)\cdot\dot{\mathbf{s}}_{i,j}(t) ~{d}t =&~ \dot{\mathbf{s}}_{i,j}(t)\cdot\dot{\mathbf{s}}_{i,j}(t) \notag \\&- \int \ddot{\mathbf{s}}_{i,j}(t)\cdot\dot{\mathbf{s}}_{i,j}(t) ~{d}t. \notag
%    \end{align}
%    Next, combine the integrals on the left hand side
%    \begin{equation}
%        2\int \ddot{\mathbf{s}}_{i,j}(t)\cdot\dot{\mathbf{s}}_{i,j}(t) ~{d}t = \dot{\mathbf{s}}_{i,j}(t)\cdot\dot{\mathbf{s}}_{i,j}(t) = a^2_{ij}(t).
%    \end{equation}
  %  Finally, take the derivative of both sides, which simplifies to
  The projection of $\ddot{\mathbf{s}}_{ij}(t)$ onto $\hat{q}_{ij}(t)$ can be calculated by applying integration by parts, which yields \cite{Beaver2020AnAgents},
    \begin{equation}
        \ddot{\mathbf{s}}_{ij}(t)\cdot\dot{\mathbf{s}}_{ij}(t) = a_{ij}(t)\cdot \dot{a}_{ij}(t). \label{eq:sddotsdot}
    \end{equation}
    Therefore, the projection of $\ddot{\mathbf{s}}_{ij}$ onto the contact basis is
    \begin{equation} \label{eq:sDdotProjection}
        \ddot{\mathbf{s}}_{ij}(t) =
        \begin{bmatrix}
        -a_{ij}^2(t)\frac{1}{2R} \\ \dot{a}_{ij}(t)
        \end{bmatrix}.
    \end{equation}
    Finally, we can use \eqref{eq:sDdotProjection} to solve for the time derivatives of \eqref{eq:phat} and \eqref{eq:qhat}. First we have
    \begin{align}
        \frac{d}{dt}\hat{p}_{ij}(t) = \frac{\dot{\mathbf{s}}_{ij}(t)}{2R} = \frac{a(t)}{2R}\hat{q}_{ij}(t). \label{eq:dPhatDt}
    \end{align}
    Taking the time derivative of $\hat{q}_{ij}(t)$ and substituting \eqref{eq:sDdotProjection} in the numerator yields \cite{Beaver2020AnAgents},
    %Then, by the quotient rule,
        \begin{align}
        \frac{d}{dt}\hat{q}_{ij}(t) %&= \frac{\ddot{\mathbf{s}}_{ij}(t)~a_{ij}(t) - \dot{\mathbf{s}}_{ij}(t)~\dot{a}_{ij}(t)}{a_{ij}^2(t)} \notag \\
        %&= \frac{a_{ij}(t)\big(-a^2(t)\frac{1}{2R}\hat{p}_{ij}(t) + \dot{a}_{ij}(t)\hat{q}_{ij}(t)\big)}{a_{ij}^2(t)} \notag\\
        %&\hspace{1em}- \frac{\dot{\mathbf{s}}_{ij}(t)~\dot{a}_{ij}(t)}{a_{ij}^2(t)} \notag \\
        %&=
        = -\frac{a(t)}{2R}\hat{p}_{ij}(t). \label{eq:dQhatDt}
    \end{align}
    
	%--- solving for the ODE
	
	By \eqref{eq:s}, we may now write $\ddot{\mathbf{s}}_{ij}(t)$ projected on to the contact basis (Definition \ref{df:contactBasis}) as
	\begin{align}
	\ddot{\mathbf{s}}_{ij}(t) =&~ \mathbf{u}_j(t) - \mathbf{u}_i(t) = -\mathbf{L}_{ij}^v(t) + m_{ij}(t)\mathbf{s}_{ij}(t) \notag\\
	=&~-\mathbf{L}_{ij}^v(t) \cdot
	\begin{bmatrix}
	\hat{p}_{ij}(t) \\ \hat{q}_{ij}(t)
	\end{bmatrix} + m_{ij}(t) \cdot
	\begin{bmatrix}
	2R \\ 0
	\end{bmatrix}, \label{eq:sddot}
	\end{align}
	where $m_{ij}(t) = \mu_{ij}(t) + \mu_{ji}(t)$, and $\mathbf{L}_{ij}^x(t) = \boldsymbol{\lambda}_j^x(t) - \boldsymbol{\lambda}_i^x(t)$. Next we substitute \eqref{eq:sddots} and \eqref{eq:sddotsdot} into \eqref{eq:sddot} and rewrite it as a system of scalar equations,
	\begin{align}
	    \mathbf{L}_{ij}^v(t)\cdot\hat{p}_{ij}(t) &= \frac{a_{ij}^2(t)}{2R} + 2R ~m_{ij}(t), \label{eq:lvp}\\
	    \mathbf{L}_{ij}^v(t)\cdot\hat{q}_{ij}(t) &= -\dot{a}_{ij}(t). \label{eq:lvq}
	\end{align}
	Taking a time derivative yields
	\begin{align}
	    \dot{\mathbf{L}}_{ij}^v(t)\cdot\hat{p}_{ij}(t) + \mathbf{L}_{ij}^v(t)\cdot\dot{\hat{p}}_{ij}(t) &= \frac{a_{ij}(t)\dot{a}_{ij}(t)}{R} + 2R \dot{m}_{ij}(t), \label{eq:lvpdot}\\
	    \dot{\mathbf{L}}_{ij}^v(t)\cdot\hat{q}_{ij}(t) + \mathbf{L}_{ij}^v(t)\cdot\dot{\hat{q}}_{ij}(t) &= -\ddot{a}_{ij}(t). \label{eq:lvqdot}
	\end{align}
	Then we substitute \eqref{eq:lambdaVdot}, \eqref{eq:dPhatDt}, and \eqref{eq:dQhatDt} into \eqref{eq:lvpdot} and \eqref{eq:lvqdot} which simplifies to
	\begin{align}
	    %\mathbf{L}_{ij}^p(t)\cdot\hat{p}_{ij}(t) &= 4R~\dot{m}_{ij}(t) - \frac{3}{2R}a_{ij}(t)\dot{a}_{ij}(t) \label{eq:lpp}\\
	    \mathbf{L}_{ij}^p(t)\cdot\hat{p}_{ij}(t) &= -2R~\dot{m}_{ij}(t) - \frac{3}{2R}a_{ij}(t)\dot{a}_{ij}(t), \label{eq:lpp}\\
	    %\mathbf{L}_{ij}^p(t)\cdot\hat{q}_{ij}(t) &= -\dot{a}_{ij}(t). \label{eq:lpq}
        \mathbf{L}_{ij}^p(t)\cdot\hat{q}_{ij}(t) &= \ddot{a}_{ij}(t) - \frac{a_{ij}^3(t)}{4R^2}. \label{eq:lpq}
	\end{align}
	%--- presenting the ODE
    %Taking an additional time derivative of \eqref{eq:lpp} and \eqref{eq:lpq} and substituting \eqref{eq:lambdaPdot}, \eqref{eq:dPhatDt}, and \eqref{eq:dQhatDt} 
    Repeating this process of deriving and substituting on \eqref{eq:lpp} and \eqref{eq:lpq} yields a pair of coupled nonlinear second-order ordinary differential equations,
    \begin{align}
    %a_{i,j}^4(t)+32R^4\ddot{m}_{i,j}(t) =&~ 12R^2\dot{a}_{i,j}^2(t) +16R^2a_{i,j}(t)\ddot{a}_{i,j}(t)\notag\\ &+ 8R^2a_{i,j}^2(t)m_{i,j}(t), \label{eq:diffeq1} \\
    %\frac{3}{2R}a_{i,j}^2(t)\dot{a}_{i,j}(t) =&~2a_{i,j}(t)\dot{m}_{i,j}(t) + \dddot{a}_{i,j}(t) %\notag\\ &+ 2m_{i,j}(t)\dot{a}_{i,j}(t). \label{eq:diffeq2}
    \frac{4a_{ij}(t)\ddot{a}_{ij}(t)}{2R} + \frac{3\dot{a}^2_{ij}(t)}{2R} &+ 2R\ddot{m}_{ij}(t) = \frac{a_{ij}^4(t)}{8R^3} \notag\\
    &~~+ \frac{m_{ij}(t)a_{ij}^2(t)}{2R}, \label{eq:diffeq1} \\
    m_{ij}(t)\dot{a}_{ij}(t) + \dot{m}_{ij}(t) a_{ij}(t) &+ \frac{6a_{ij}^2(t)\dot{a}_{ij}(t)}{4R^2} = \dddot{a}_{ij}(t). \label{eq:diffeq2}
    \end{align}
    Thus, for any constrained trajectory to be energy-optimal it must be a solution of \eqref{eq:diffeq1} and \eqref{eq:diffeq2} while also satisfying the boundary conditions \eqref{eq:IC} and \eqref{eq:phiBoundary}.
    In general this is difficult, as both equations are nonlinear and \eqref{eq:diffeq2} is third order.
    
    An alternative solution is to impose $a_{ij}(t) = 0$ over any nonzero interval where the safety constraint is active. This implies $\dot{\mathbf{s}}_{ij}(t) = \mathbf{v}_j(t) - \mathbf{v}_i(t) = 0$. Additionally, as $\dot{\mathbf{s}}_{ij}(t)$ is constant, its derivative $\ddot{\mathbf{s}}_{ij}(t) = 0$. Thus we may select $\mathbf{v}_i(t) = \mathbf{v}_j(t)$ and $\mathbf{u}_i(t) = \mathbf{u}_j(t)$ as a ``reigning optimal'' solution \cite{Ross2015}.

	To generate a final energy-optimal trajectory, each agent will piece together unconstrained and constrained arcs during each trajectory update. In the following section, we present the decentralized strategy used by each agent to generate collision-free by piecing together energy-optimal motion primitives.
	
	\subsection{Decentralized Trajectory Generation} \label{ss:generation}
	
	To generate the trajectory of each agent we will use the following definition of contact intervals.
	\begin{definition} \label{df:contactIntervals}
		For each agent $i\in\mathcal{A}$ we define $g\in\mathbb{N}$ \emph{contact intervals} indexed by $k = 1, 2, ..., g$. These $g$ intervals must satisfy
		\begin{align}
		\tau_i^k &\subset [t_0, t_f] \text{ such that } \notag\\
		\mathcal{V}_i(t_a) &= \mathcal{V}_i(t_b) \neq \emptyset ~~ \forall~ t_a, t_b \in \tau_i^k, \notag
		\end{align}
		where for any two intervals $p,q\in\mathbb{N}, ~ p\neq q$, we have $t_p < t_q$ for $t_p\in\tau_i^p$ and $t_q\in\tau_i^q$. The contact intervals correspond to the nonzero and non-overlapping intervals of time where $\mathcal{V}_i(t)$ is invariant.
	\end{definition}
	The contact intervals correspond to the instances in time when each agent $i\in\mathcal{A}$ would violate the safety constraint when traveling along an unconstrained trajectory. To guarantee safety, we enforce the centralized solution to the safety-constrained case over these intervals. We then piece these constraints together with the initial and final conditions using unconstrained trajectory segments. Each agent $i\in\mathcal{A}$ performs the following steps simultaneously:
	\begin{enumerate}
		\item Generate an unconstrained trajectory for $t\in[t_0, t_f]$.
		\item Exchange unconstrained trajectories with all $j\in\mathcal{N}_i(t)$.
		\item Calculate $\mathcal{V}_i(t)$.
		\item Generate a constrained arc for every contact interval (Definition \ref{df:contactIntervals}) and fix the entry and exit states for each interval.
		\item Piece together all constrained arcs using unconstrained trajectories and continuity in $\mathbf{x}_i$ at the boundaries.
		\item Generate initial and final unconstrained arcs to satisfy $\mathbf{x}_i(t_0)$ and $\mathbf{x}_i(t_f)$.
		\item Generate escape arcs, i.e., agent $i$ may exit to an unconstrained trajectory early if it does not violate any safety constraints.
	\end{enumerate}
	The above steps will be performed by all agents simultaneously and repeated with a period of $\Delta T$ per our control framework. %For any agent $i\in\mathcal{A}$, the trajectory of all $j\in\mathcal{N}_i(t)$ are estimated from their most recent trajectory. % The details of each step are presented for some agent $i\in\mathcal{A}$ next, but all steps are performed by all agents simultaneously.

	\section{Simulation Results} \label{sc:simulation}
	
	To validate our decentralized controller we developed two sets of simulations. First, we implemented the controller in MATLAB where the unconstrained and distance-constrained trajectories could be validated on double integrator agents.
	Next, the controller was applied to a set of AscTec quadrotors in Gazebo. The results of this simulation show that our optimal controller generates high-level trajectories which display emergent flocking behavior. These trajectories are also realizable by the dynamics and low-level flight controller of a commercially available quadrotor.
	
	%To enhance the performance of the controller, the coefficients in the flocking error function, \eqref{eq:flockError}, were normalized such that each component \eqref{eq:velocityControl} - \eqref{eq:aggregation} had a magnitude around 1 during normal operation. These scaling factors are given by
	%\begin{align}
%		w_1 = \frac{k_d}{||\mathbf{v}_d||}, ~~~
%	    w_3 = \frac{k_p}{D^2 N}, ~~~
%	    w_2 =  \frac{k_v}{||\mathbf{v}_d||},
%	\end{align}
%	where $k_p, k_v, $ and $k_d$ are control gains, and are given along with the remaining system parameters in Table \ref{tab:parameters}.
	
%    \begin{table}[ht]
%        \centering
%        \vspace{1em}
%        \begin{tabular}{l|cc}
%            Variable & Symbol & Value \\ \toprule
%            Horizon & $T$ & 0.4 s \\
%            Update Period & $\Delta T$ & 0.1 s\\
%            Position gain & $kp$ & 1.5 \\
%            Velocity gain 1 & $kv$ & 1.0 \\
%            Velocity gain 2 & $kd$ & 3.0
%        \end{tabular}
%        \caption{System parameters used for the MATLAB and Gazebo simulations.}
%        \label{tab:parameters}
%    \end{table}
    
\begin{figure}[ht]
		\centering
		\includegraphics[width=0.49\linewidth]{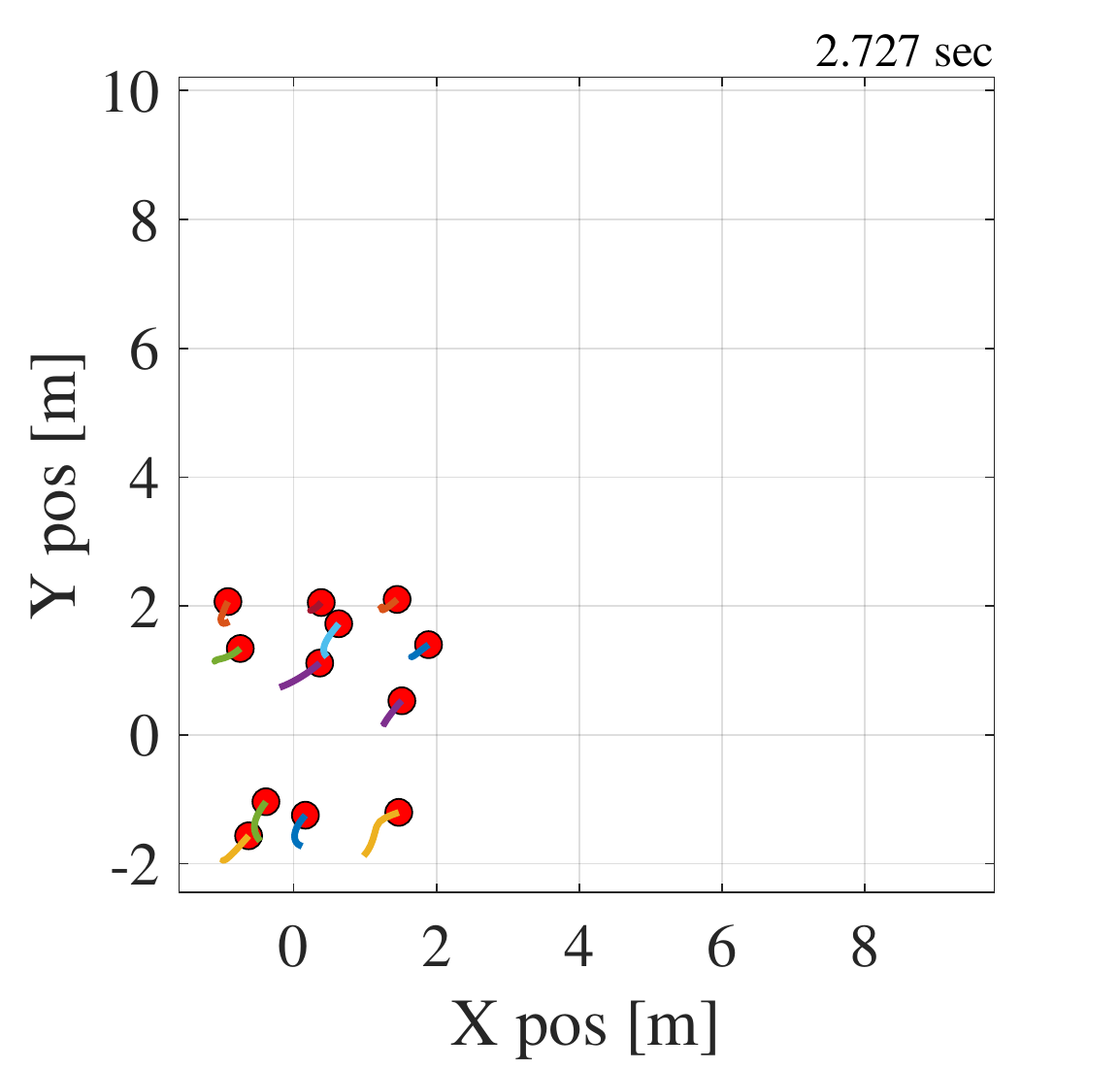}
		\includegraphics[width=0.49\linewidth]{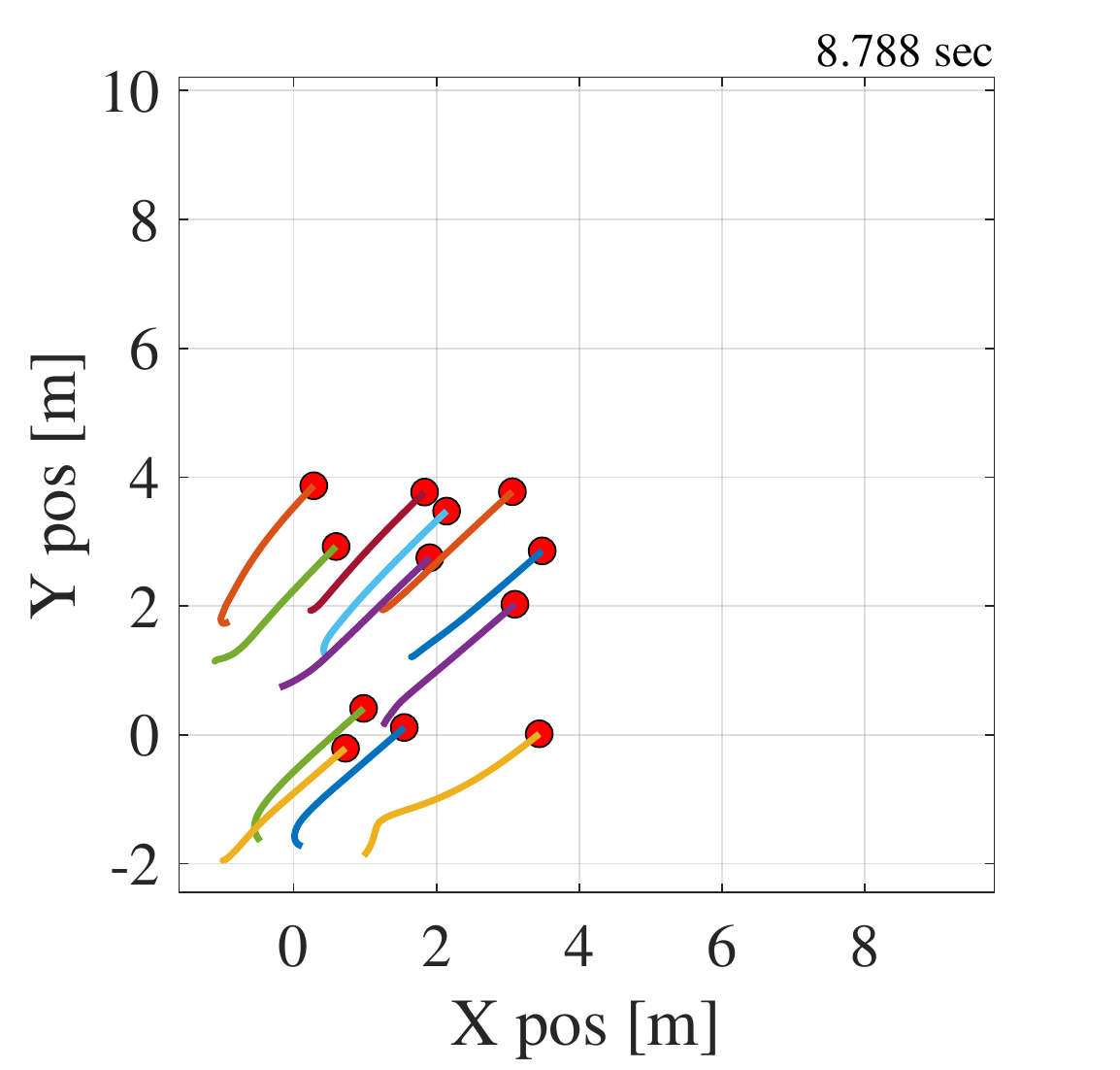}
		\includegraphics[width=0.49\linewidth]{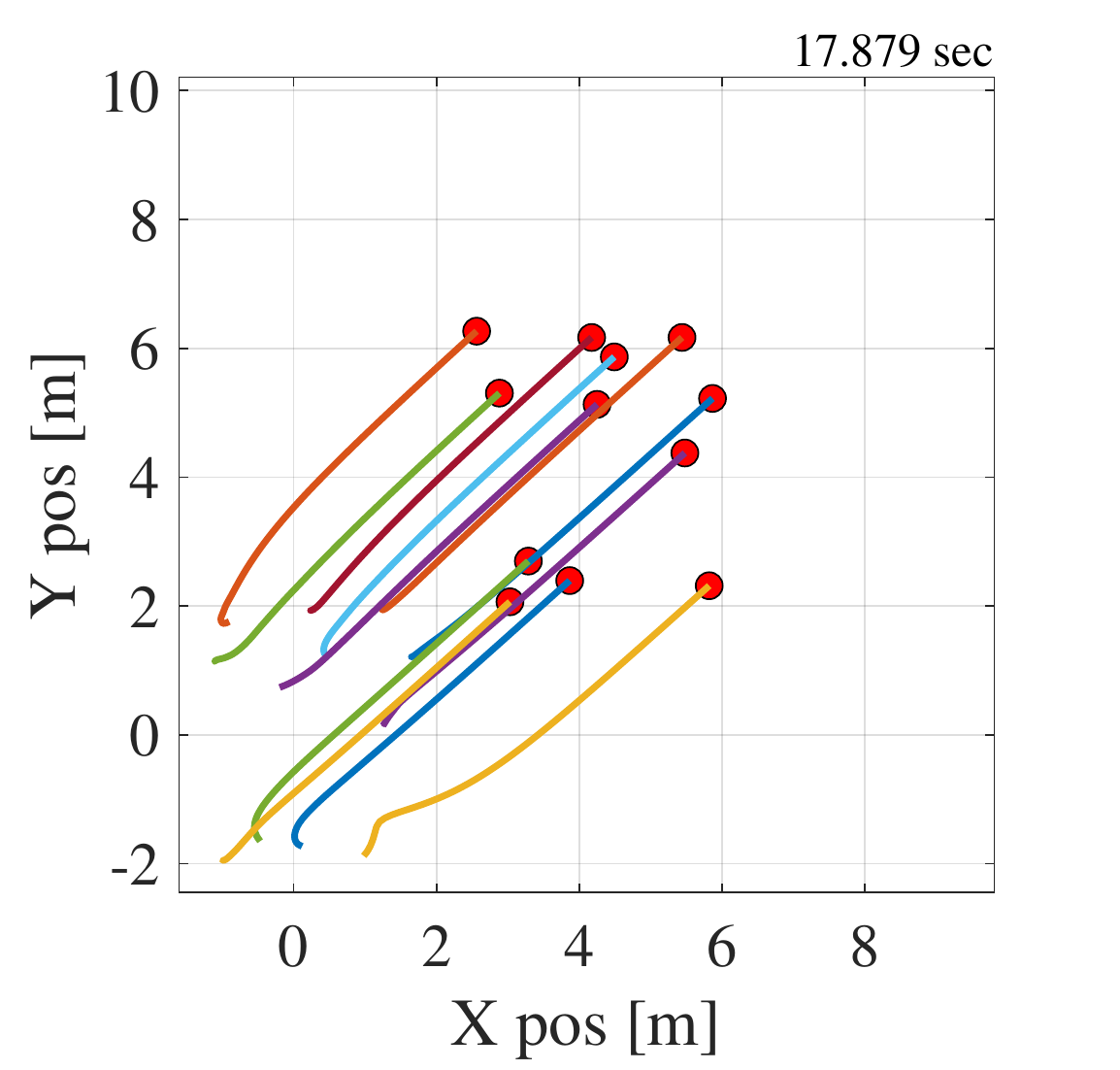}
		\includegraphics[width=0.49\linewidth]{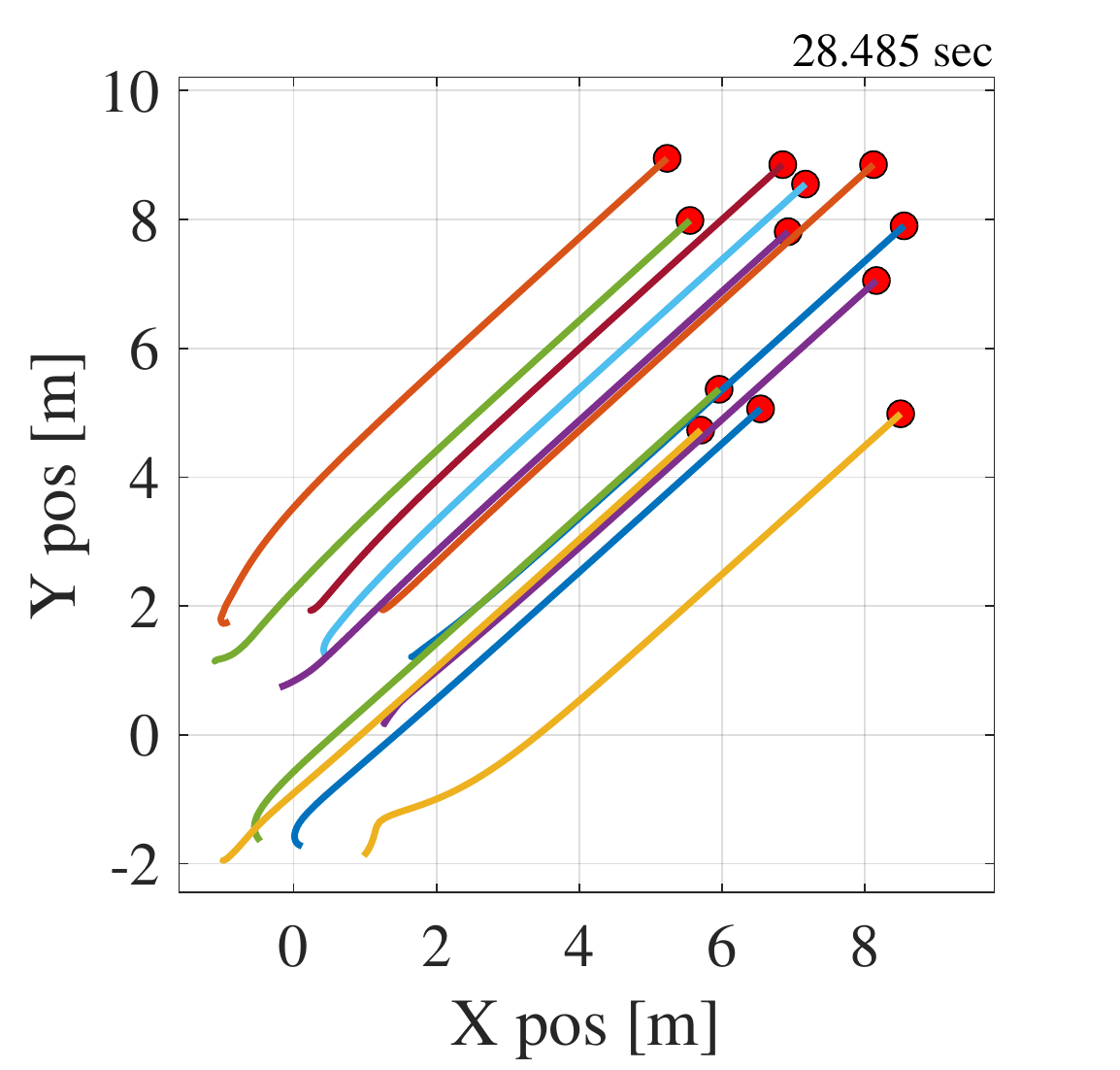}
		\caption{Flocking behavior for 12 agents simulated in MATLAB.}  \label{fg:matlabFlock}
	\end{figure}

%	\subsection{MATLAB Simulation}
   % We explored two scenarios through simulation in MATLAB and Gazebo. %We first initialized two agents such that a collision would occur in the unconstrained case  (Figure \ref{fg:matlabCollision}). This result illustrates the collision-avoidance capability of the proposed algorithm. Next 
    First, we validated our proposed controller by placing $12$ agents at feasible initial points randomly within the domain in MATLAB. The resultant flocking behavior if presented in Fig. \ref{fg:matlabFlock}; this shows that after a short transient period, a stable flock is formed which moves to the northeast as specified by $\mathbf{v}_d$.
    
    %and \ref{fg:matlabStates}. %Figure \ref{fg:matlabStates} shows an initial spike in energy consumption during the transient flock formation. The energy consumption then quickly approaches zero as a stable flocking aggregate is formed.
	
%	\begin{figure}[ht]
%		\centering
%		\includegraphics[width=0.9\linewidth]{matlabUvT.pdf}
%		\caption{A plot of maximum and average acceleration magnitude vs time for the flocking behavior in MATLAB (Figure \ref{fg:matlabFlock}).}  \label{fg:matlabStates}
%	\end{figure}

%    The results presented in Figures \ref{fg:matlabFlock} and \ref{fg:matlabStates} show that after a short transient period a stable flock is formed. This flock moves in the direction specified by $\mathbf{v}_d$ and begins to cruise at a constant speed around $t=10$ seconds. This smooth, stable behavior results from each agent following a minimum energy trajectory while forming the flock aggregate.

\begin{figure}[ht]
		\centering
		\includegraphics[width=0.49\linewidth]{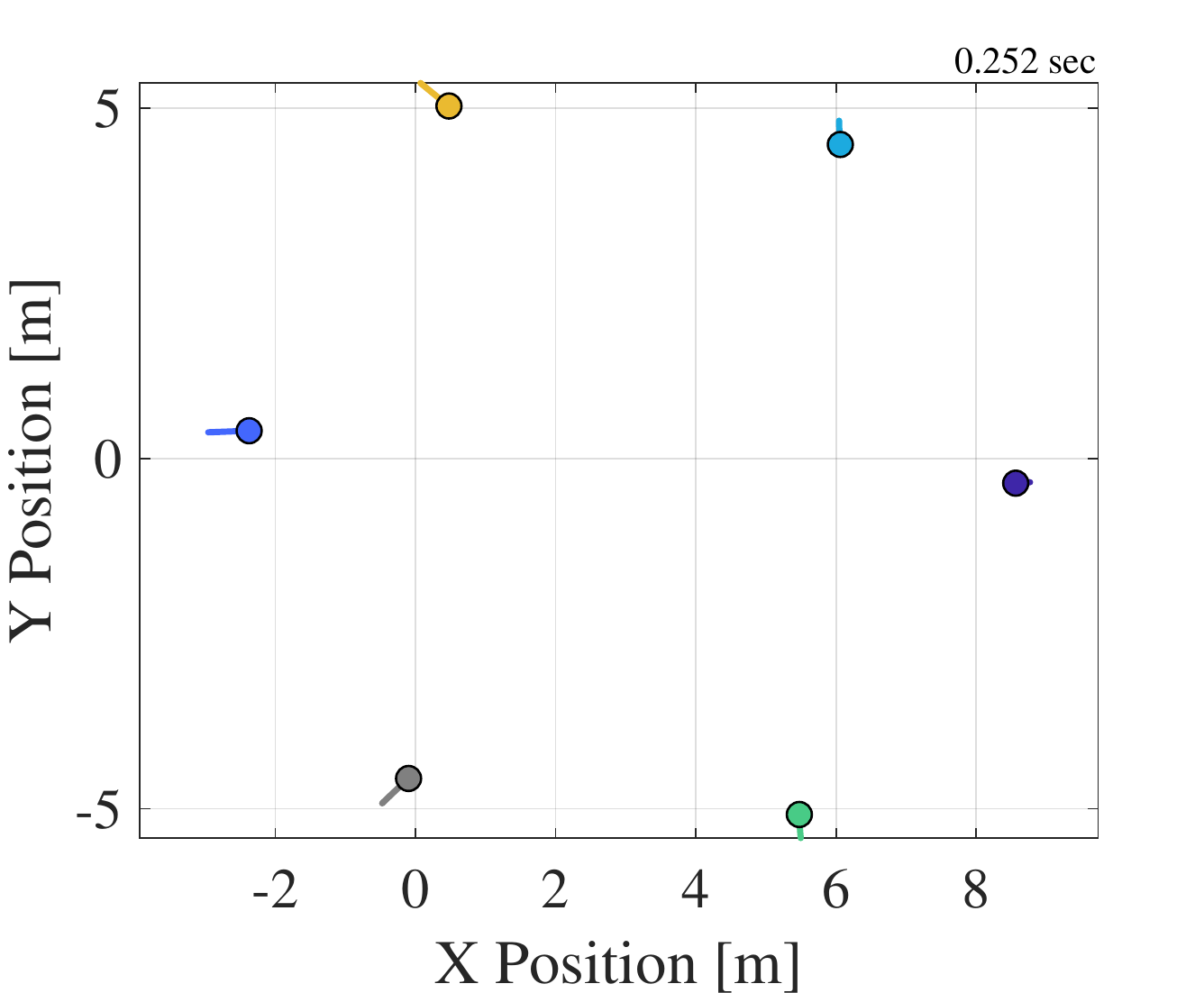}
		\includegraphics[width=0.49\linewidth]{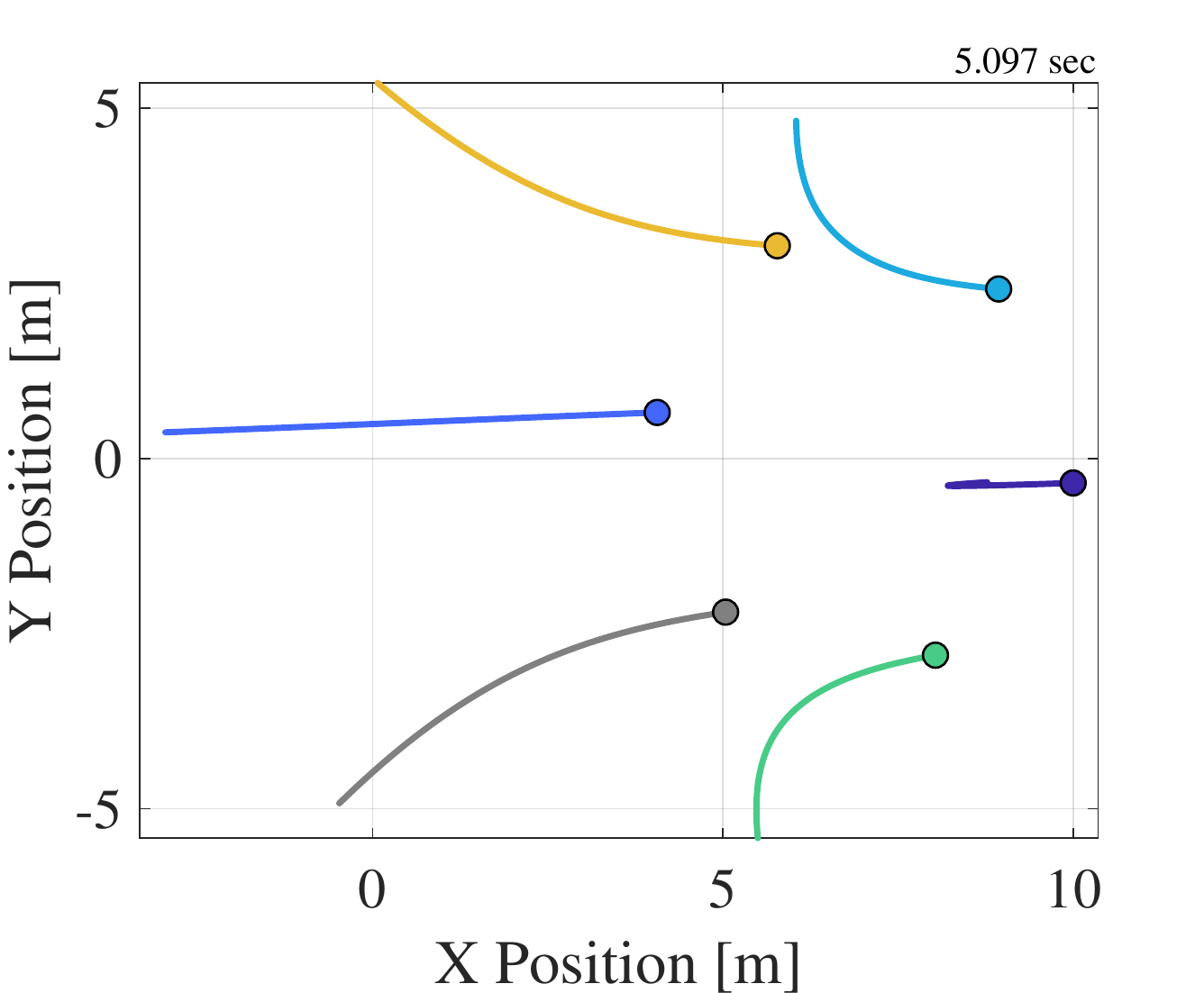}
		\includegraphics[width=0.49\linewidth]{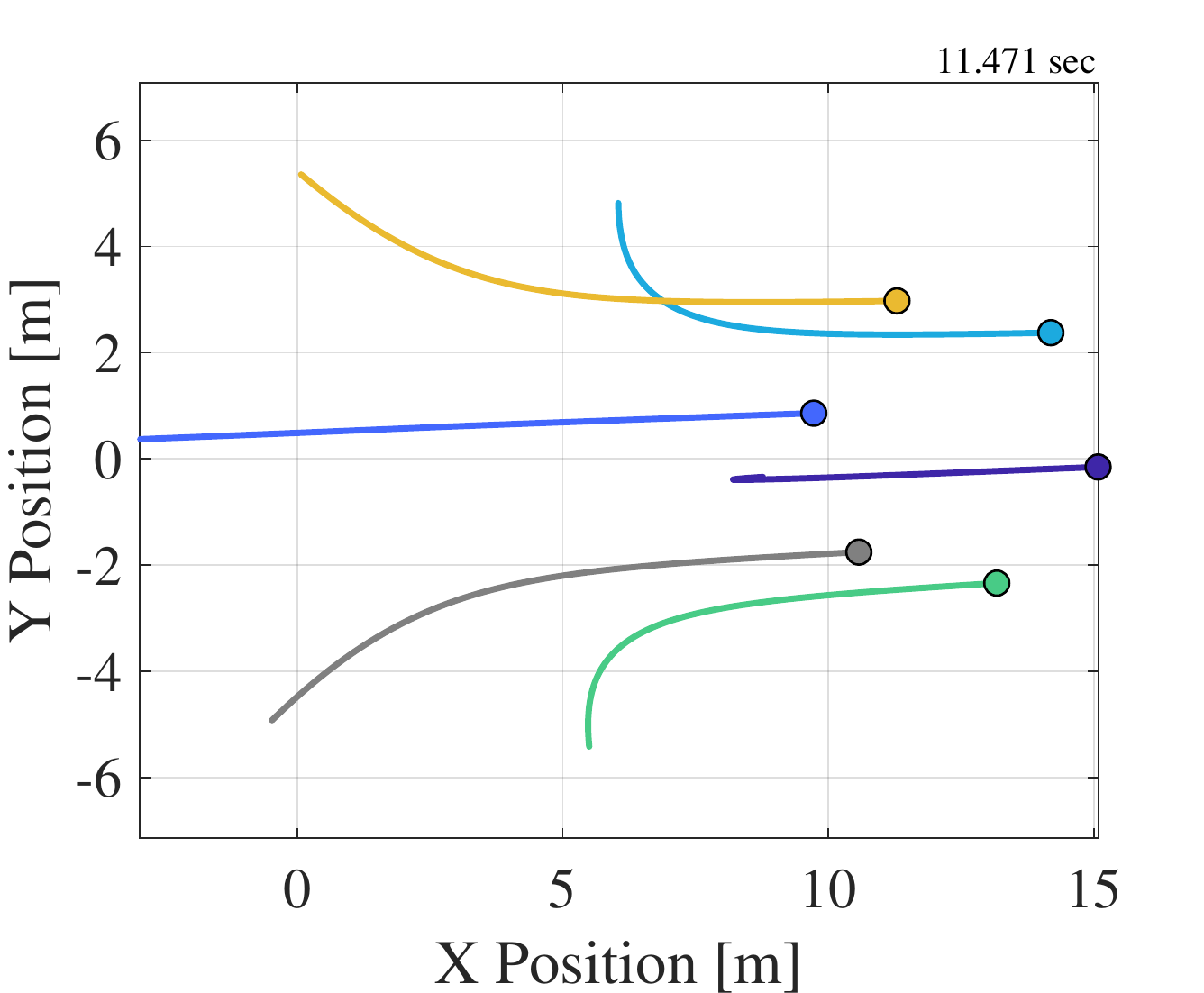}
		\includegraphics[width=0.49\linewidth]{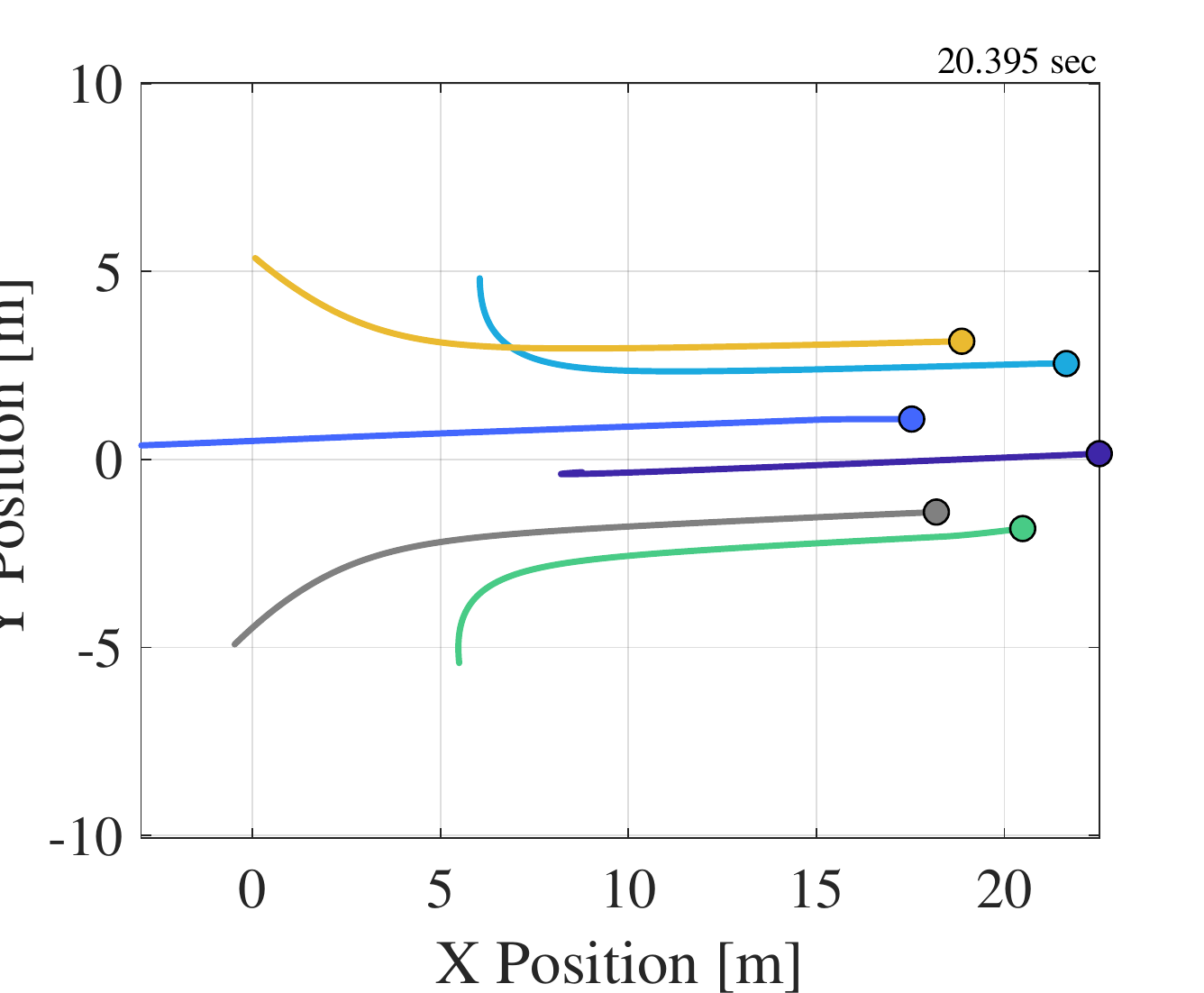}
		\caption{Gazebo simulation with six agents spaced in a hexagonal formation.}  \label{fg:gazeboFlock}
	\end{figure}

\begin{figure}[ht]
		\centering
		\includegraphics[width=0.95\linewidth]{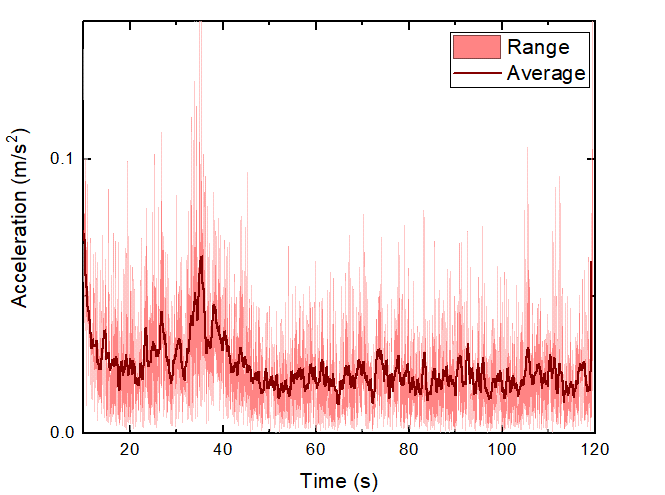}
		\caption{Time series of the average and maximum acceleration magnitude for the six drones.}  \label{fg:gazeboEnergy}
	\end{figure}

%	\subsection{Gazebo Simulation}
	Next, we implemented a \emph{sensing-based approximation} of the derived optimal control algorithm in Gazebo. This approximation does sacrifice optimality; however, it does not require communication between the agents and it allows an explicit closed-form solution of the boundary conditions. The sensing-based approach for agent $i\in\mathcal{A}$ approximates \eqref{eq:boundary3} as
	\begin{align}
	    \mathbf{p}_i(t_f) =&~ \sum_{j\in\mathcal{N}_i}\mathbf{p}_{j}(t_0) + \frac{\mathbf{a}_i}{2w_3\sum_{j}(||\mathbf{s}_{ij}(t_0)|| - D)}, \label{eq:boundary3-s}
	\end{align}
	where $\mathbf{p}_j(t_0), ~ j\in\mathcal{N}_i(t_0)$ must only be sensed by agent $i$. The expected flocking behavior still emerges and is presented in Fig. \ref{fg:gazeboFlock}.
    We implemented our controller in Gazebo using the RotorS package \cite{Furrer2016RotorS---AFramework} and six AscTec Hummingbirds operating in a horizontal plane. The parameters used for this simulation were: $D = 0.25$ m, $\mathbf{v}_d = (2.5, 0)$ m/s, and $h = 4.5$ m. The results of these simulations are visualized in Figs. \ref{fg:gazeboFlock} and \ref{fg:gazeboEnergy}.
     In Fig. \ref{fg:gazeboFlock}, agents initially coalesce into a hexagonal pattern before flocking along the direction of $\mathbf{v}_d$ ($+x$). The flock naturally forms a hexagonal near the end of the simulation due to the relative distance term in the cost function. Figure \ref{fg:gazeboEnergy} shows the energy consumption of the agents decaying despite the use of the sensing approximation to estimate neighbor trajectories.

    The Gazebo simulation shows an asymptotic reduction in energy consumption for agents with significantly more complicated dynamics, i.e., drones. Additionally, the sensing-only approximation \eqref{eq:boundary3-s} does not seem to have a significant impact on the global structure of the flock and rate of energy consumption in this scenario.

    \section{Conclusion} \label{sc:conclusion}
    
    In this paper, we proposed an optimal control approach to realize flocking behavior in a group of cooperative agents. We presented the optimal trajectory in the form of a boundary value problem and provided a candidate solution which is locally optimal. 
    Some potential directions for future research include extending the obstacle avoidance constraint to include terrain and mobile obstacles, redefining the aggregation function to use a local description of flocking, and testing the optimal control algorithm on physical hardware in the University of Delaware's Scaled Smart City.
	
	\section*{Acknowledgments}
	The authors would like to thank Michael Dorothy at Combat Capabilities Development Command, Army Research Laboratory, for his insightful discussions on optimal control.
	
    \bibliographystyle{IEEEtran}
    \bibliography{mendeley,Andreas,IDS_Publications_030920}

\end{document}